\documentclass[final,11pt,times]{elsarticle}

\usepackage{graphicx}
\usepackage[font=bf,labelfont=bf]{caption}
\usepackage{amssymb}
\usepackage{amsmath}
\usepackage{multirow}
\usepackage{algorithm}
\usepackage{algpseudocode}
\usepackage{url} 
\usepackage{placeins}
\usepackage{geometry}
\usepackage{soul}
\usepackage{xcolor}
\usepackage{float}
\usepackage{fancyhdr}
\usepackage[normalem]{ulem}

\journal{CIE51 Proceedings}

\begin{document}

\begin{frontmatter}



\title{A Bi-Objective Mathematical Model for the Multi-Skilled Resource-Constrained Project Scheduling Problem Considering Reliability: An AUGMECON2VIKOR Hybrid Method}


\author[inst1]{Mohammad Ghasemi}

\affiliation[inst1]{organization={School of Information Technology},
            addressline={Deakin University}, 
            city={Geelong},
            postcode={3216}, 
            state={VIC},
            country={Australia}}

\author[inst1]{Asef Nazari}
\author[inst1]{Dhananjay Thiruvady}
\author[inst2]{Reza Tavakkoli-Moghaddam}
\author[inst2]{Reza Shahabi-Shahmiri}
\author[inst3]{Seyed-Ali Mirnezami}

\affiliation[inst2]{organization={School of Industrial Engineering, College of Engineering}, 
            addressline={University of Tehran},
            city={Tehran},
            country={Iran}}

\affiliation[inst3]{organization={Department of Industrial Engineering, Faculty of Engineering},
            addressline= {Shahed University},
            city={Tehran},
            country={Iran}} 

\begin{abstract}
In recent years, resources with multiple skills have received attention as an extension of the resource-constrained project scheduling problem known as MSRCPSP. Although the disruption rate is well-estimated in today's manufacturing projects, its impact on project makespan and cost need further investigation. Hence, this study presents a novel mathematical model for the MSRCPSP considering reliability, namely MSRCPSPR. The model proposes both objectives of minimizing project makespan and project cost.
The MSRCPSP is an NP-hard problem, and including reliability constraints, as proposed in this paper, makes solving the problem more intractable. 
To cope with the computational challenges of solving the problem, a combination of an enhanced version of the epsilon-constraint method as well as an augmented version of the VIKOR algorithm,
namely AUGMECON2VIKOR, is employed to solve benchmark instances j10 and j20 from the PSPLIB. A comparative analysis demonstrates the performance of the proposed method, and the sensitivity analysis represents the effects of positive reliable constraints on the objective functions. Employing the proposed method, the project makespan and cost are reduced by nearly 2.55\% and 2.80\% in j10 on average. CPU time is also decreased by about 543 seconds in comparison to the epsilon-constraint method.
\end{abstract}


\end{frontmatter}

\thispagestyle{fancy}


\section{INTRODUCTION}
\label{sec:introduction}
The resource-constrained project scheduling problem (RCPSP) has been one of the most common issues in operation management for many decades \citep{ghasemi2023chance}. In recent works, the objectives have mainly been to optimize makespan and financial aspects \citep{thiruvady2014lagrangian, thiruvady2019maximising}. Among various extensions, the multi-skilled version (MSRCPSP) represents a more realistic representation of real-world projects. This problem was first addressed by \citet{drezet2008project}. In the MSRCPSP, activities need specific skills to be carried out, and multi-skilled resources have these necessary skills. This problem has practical applications in a wide range of domains, such as research and development \citep{wang2022multi}, consulting \citep{zarei2024multi}, construction \citep{joshi2019effective}, IT \citep{chen2020competence, heimerl2010scheduling}, and manufacturing \citep{chen2022multi}. The numerous publications in this context \citep{afshar2021multi} demonstrate its relevance and importance. 

To the best of our knowledge, most project scheduling developments have been proposed without considering any risk of disruption of renewable resources. The risk of failure in machine components is the most vital factor that has an impact on project scheduling and total project cost \citep{souier2019nsga}. To justify the initial investment, project managers need to gain a high level of reliability and risk management considerations to cope with unforeseen circumstances. Hence, to deal with the aforementioned issues in the MSRCPSP, especially in manufacturing projects, this paper presents a novel bi-objective mathematical model for project scheduling to minimize the project makespan and cost simultaneously with reliability consideration, namely MSRCPSPR. This model also takes into account the waiting time for renewable resources using queuing theory concepts.


The MSRCPSP is vastly investigated in the literature, and here, a short review of the most recent works is provided. \citet{zarei2024multi} studied multi-skilled workforce and their skill levels to maximize quality and minimize project makespan in addition to resource and tardiness costs under uncertain conditions. Also, \citet{torba2024solving} proposed a time-indexed mixed-integer linear programming model for the multi-project MSRCPSP to minimize the weighted tardiness of projects and the weighted duration of projects. Learning of human resources in the MSRCPSP was investigated by \citet{mozhdehi2024multi} to find the best start and finish time of project activities utilizing a modified discrete variant of the biogeography-based optimization (MBBO) algorithm. Benchmark instances for the MSRCPSP were addressed by \citet{snauwaert2023classification}. The instances were validated through software and railway construction companies. In addition to the project makespan and project costs, minimizing project risk with stochastic non-renewable resources was addressed by \citet{mirnezami2023integrated} in MSRCPSP using an AUGMECON2 method. However, in these studies, the main focus was on workforce skills in the project resources, while in many real-world projects, machines are utilized to perform activities. Thus, considering machines as renewable resources with multiple skills is a crucial issue for projects in numerous industries. Including machines in addition to the workforce as renewable resources in MSRCPSP provides possible conditions to create queues due to the disruption. 

In the mentioned articles, researchers have only focused on workforce skills in the project resources, while in many projects, machines are utilized to perform activities. Thus, considering machines as renewable resources with multiple skills is a crucial issue for projects in various industries. Including machines in addition to the workforce as renewable resources in the MSRCPSP provides possible conditions to create queues due to the disruption. Hence, calculating the waiting time of machines within the project makespan is substantial, which has not been investigated in project scheduling problems. To address the gaps in recent literature, the novelties of this work are summarised as follows:
\begin{itemize}
    \item Present a novel MINLP model to include the MSRCPSP considering reliability constraints.
    \item Calculate waiting time for each renewable resource in the MSRCPSP for the first time. 
    \item Implement an AUGMECON2VIKOR algorithm to solve the bi-objective model for the MSRCPSPR.
    \item Consider arrival and retrieval rates based on the $M/M/1$ queue system in the project scheduling problem.
\end{itemize}
    

\section{PROBLEM DEFINITION}
\label{sec:problem definition}
This section describes the definition of the MSRCPSP in addition to the waiting model when an unexpected disruption leads to waiting time in the project activities. The MSRCPSP considers both scheduling activities and allocating resources. This problem can be considered as a combination of two separate sub-problems as outlined in \citep{zheng2017teaching,wang2018knowledge}. In general, there are two common representation networks to demonstrate precedence relationships between activities in a project, which have been previously investigated. In this study, the activity-on-node (AON) network is selected due to a clearer representation of project activity relationships as emphasized by \citet{kolisch2001integrated}. A formal definition of each subproblem is as follows.

This problem considers a single project with a topologically ordered acyclic AON network, in which $N=\{1,\ldots,|N|\}$ is the set of nodes that represent the project activities. The start and end nodes are indicated by dummy activities. Therefore, there are $|N|-2$ activities that can be executed and denoted by $I=N-\{1, |N|\}$. The processing time of dummy nodes is zero and no resource or skill is required for them. The processing time of each project activity $i$ is indicated by $d_i$. The time lag between precedence activities is considered to be zero. Furthermore, activities cannot be preempted or crashed. For project execution, activities needs multi-skilled renewable resource $k$ from the resource set $R=\{1,\ldots,|R|\}$ and skill $s$ from the skill set $S=\{1,\ldots,|S|\}$. In this study, equal efficiency for all members of the resource set is assumed.  Each resource $k$ can master a range of skills. $b_{lk}$ is defined to denote whether resource $k \in R$ masters skill $l \in S$. It should be mentioned that each activity $i \in I$ requires $r_{il}$ resources to be performed. Each resource $k \in R$ cannot be allocated to more than one activity at each time. All the parameters are assumed to be deterministic values.

Consider multi-skilled machines as renewable resources. As there are numerous components in the execution of a project, a negligible amount of uncertainty may lead to tremendous disruption \citep{song2024impact}. Consequently, significant costs are incurred in the completion of the project. Uncertainties like electricity disruptions, unpredicted work overloads, and human resources strikes can affect expenses and completion time. Thus, disruptions and unavailability of resources are probable in the life span of a project. A decrease in the retrieval rate of resources in these conditions gives rise to extra expense and time. 

Many researchers have studied the use of queuing approach to appraise the waiting time of operation. For instance, \citet{souier2019nsga, souza2020survey, may2021queue, yu2024production} suitably employed queuing approaches to model the waiting time of product flows. However, to the best of our knowledge, no study considers the queuing theory in the project scheduling problem. Based on the real-world project in these studies, the Poisson distribution has been taken into account as a suitable model of arrival rates, in which there is a variation around scheduled arrival times \citep{peterson1995models}. This distribution can be considered in other service systems, where arrival rates, as well as capacity levels, vary remarkably over time. 


To construct a disruption model compatible with the real-world circumstances in the MSRCPSP, stochastic disruption of the queue system at each machine is considered, and it is assumed that the disruption is retrieved with special rates. Employing the $M/M/c$ queue system, independent service times are considered for each resource $k \in R$ which are distributed exponentially with rate $\mu_k$. During disruptions, the number of operational servers decreases from $c$ to $c'$ and the service rates alter from $\mu_k$ to $\mu_k'$. After retrieving the machine, the number of working servers as well as the related service rates are restored to $c$ and $\mu_k$. It is assumed that disruptions are according to a Poisson distribution with rate $\upsilon_k$, and the retrieve times are i.i.d. exponential with rate $r_k$. The flow arrivals follow a homogeneous Poisson distribution with intensity $\lambda$. Based on the abovementioned assumptions, the $M/M/1$ queue system is taken into account for all machines in this investigation.

\section{MATHEMATICAL FORMULATION}
\label{sec:mathematical formulation}
In this section, an MINLP model is presented based on \citet{snauwaert2023classification} to represent the MSRCPSPR with the disruption considerations. Including a waiting time setting based on an $M/M/1$ queue system makes the model realistic and compatible with complex real-world projects.

\subsection{Notations}
\noindent Sets:\\[0.3cm]
\begin{tabular}{ll}
    $I$ & Executable activities \\
    $N$ & Project activities including dummy activities \\
    $S$ & Skills \\
    $R$ & Resources
\end{tabular}\\[0.3cm]
Parameters:\\[0.3cm] 
\begin{tabular}{ll}
    $r_{il}$ & The number of required resources for activity $i$ with skill $l$ \\
    $b_{lk}$ & $1$ if resource $k$ has skill $l$; otherwise $0$ \\
    $p_{ij}$ & $1$ if activity $i$ precedes activity $j$; otherwise $0$\\ 
    $d_i$ & Processing time of activity $i$ \\
    $r_k$ & Retrieval time rate of resource $k$ \\
    $\upsilon_k$ & The disruption rate of resource $k$.\\
\end{tabular}\\[0.3cm]
\begin{tabular}{ll}
    $\mu_k$ & Service time at resource $k$ \\
    $c_{lk}$ & Cost of skill $l$ in using resource $k$ \\
    $M$ & A big number
\end{tabular}\\[0.3cm]
Binary Decision Variables:\\[0.3cm]
\begin{tabular}{ll}
    $X_{ilk}$ & $1$ if activity $i$ with skill $l$ performed by resource $k$; otherwise 0 \\
    $Z_{ij}$ & $1$ if activity $i$ precedes activity $j$; otherwise 0 \\
    $Y_{ik}$ & $1$ if activity $i$ requires resource $k$; otherwise 0
\end{tabular}\\[0.3cm]
Positive Decision Variables:\\[0.3cm]
\begin{tabular}{ll}
    $\lambda_k$ & Arrival rate of activities to each resource $k$ \\
    $W_k$ & Waiting time of each resource $k$\\
    $S_i$ & Start time of activity $i$\\
    $T_{i}$ & Maximum waiting time of activity $i$ 
\end{tabular}
\subsection{The Mathematical Model of the MSRCPSPR}
\vspace{-1.5em}
\begin{align}
 &\hspace{-0.9cm} \mathbf{P_1:} \;\;\hspace{0.07cm} Z_{1} = \min\ S_{|N|} \label{obj1}  \\
 & Z_{2} = \min\ \sum_{i \in I}\sum_{l \in S}\sum_{k \in R} d_i\  c_{kl}X_{ilk} \label{obj2}\\
& \text{s.t.} \notag\\
&\sum_{k\in R}X_{ilk} = r_{il} \hspace{6cm} \forall\ {i \in I,l \in S} \label{rest1} \\
&\sum_{l\in S}X_{ilk} \le 1 \hspace{6.1cm} \forall\ {i\in I,k \in R} \label{rest2}\\
&Z_{ij} + Z_{ji} \le 1 \hspace{6cm} \forall\ i,j \in N ;\ i < j \label{rest3}\\
&\sum_{l\in S}X_{ilk} + \sum_{l\in S}X_{jlk} \le 1+Z_{ij} + Z_{ji} \hspace{2.85cm}  \forall\ i,j \in N ;\ i < j, \ k \in R \label{rest4} \\
&\sum_{l \in S}\sum_{i \in I}X_{ilk} = \lambda_k \hspace{5.4cm} \forall\ k \in R \label{rest5}\\
&W_k = \frac{(r_k+\upsilon_k)^2+\mu_k\upsilon_k}{(r_k+\upsilon_k)(r_k\mu_k-r_k\lambda_k-\lambda_k\upsilon_k)}  \hspace{2.55cm} \forall\ k \in R \label{rest6}\\
&x_{ilk} \le M \cdot Y_{ik} \hspace{5.95cm} \forall\ i \in I,l \in S,k \in R \label{rest7}\\
&W_k \cdot Y_{ik}\le T_i \hspace{6cm} \forall\ {i \in I,k \in R} \label{rest8} \\
&S_i+d_i+T_i-M(1-Z_{ij})(1-p_{ij}) \le S_{j} \hspace{2cm} \forall\ i,j \in N \label{rest9}\\
&X_{ilk} \le b_{lk} \hspace{6.6cm} \forall\ {i \in N,l \in S,k \in R} \label{rest10} \\
&X_{ilk},\ Z_{ij},\ Y_{ik} \in \{0,1\};\quad   \lambda_k,\ W_k,\ S_i,\ T_i \ge 0 \hspace{1.4cm} \forall\ i,j \in N;\ \forall\ {l \in S,k \in R}  \label{rest11}
\end{align}




The first objective function $Z_1$ in Eq. (\ref{obj1}) aims to minimize the project makespan. Total project cost is minimized in the second objective function $Z_2$ in Eq. (\ref{obj2}). The skill requirements of each project activity are described by Eq.(\ref{rest1}). In addition, Eq. (\ref{rest2}) ensures that each resource $k$ can execute at most one skill in each activity. The logical precedence relationship between activities is described by Eq.(\ref{rest3}). Furthermore, Eq. (\ref{rest4}) assures the disjunctive resource allocation. Once two project activities $i$ and $j$ are scheduled to be performed in parallel, the right-hand side of this constraint is equal to $1$, as a result, $\sum_{l\in S}X_{ilk} + \sum_{l\in S}X_{jlk}$ needs to be smaller than or equal to $1$. Accordingly, the resource $k$ can merely perform one skill on activity $i$ or $j$ and will not execute two activities simultaneously. The arrival rate of each project activity is determined by Eq. (\ref{rest5}). Waiting time for each machine as a renewable resource is addressed by Eq. (\ref{rest6}), and Eq. (\ref{rest7}) ensures that each activity can be performed by resource $k$ with skill $l$ only when resource $k$ has been allocated. The waiting time for each renewable resource is limited by Eq. (\ref{rest8}), and Eq. (\ref{rest9}) guarantees the minimum finish-to-start precedence relations with a zero time lag. Finally, Eq. (\ref{rest10}) represents that a resource $k$ can merely execute a skill $l$ that it has mastered, and Eq. (\ref{rest11}) denotes the domain of decision variables.

\section{SOLUTION STRATEGY}
\label{sec:solution strategy}    
\sethlcolor{lightgray}

\begin{figure}[h]
    \centering
    \includegraphics[width=1\linewidth]{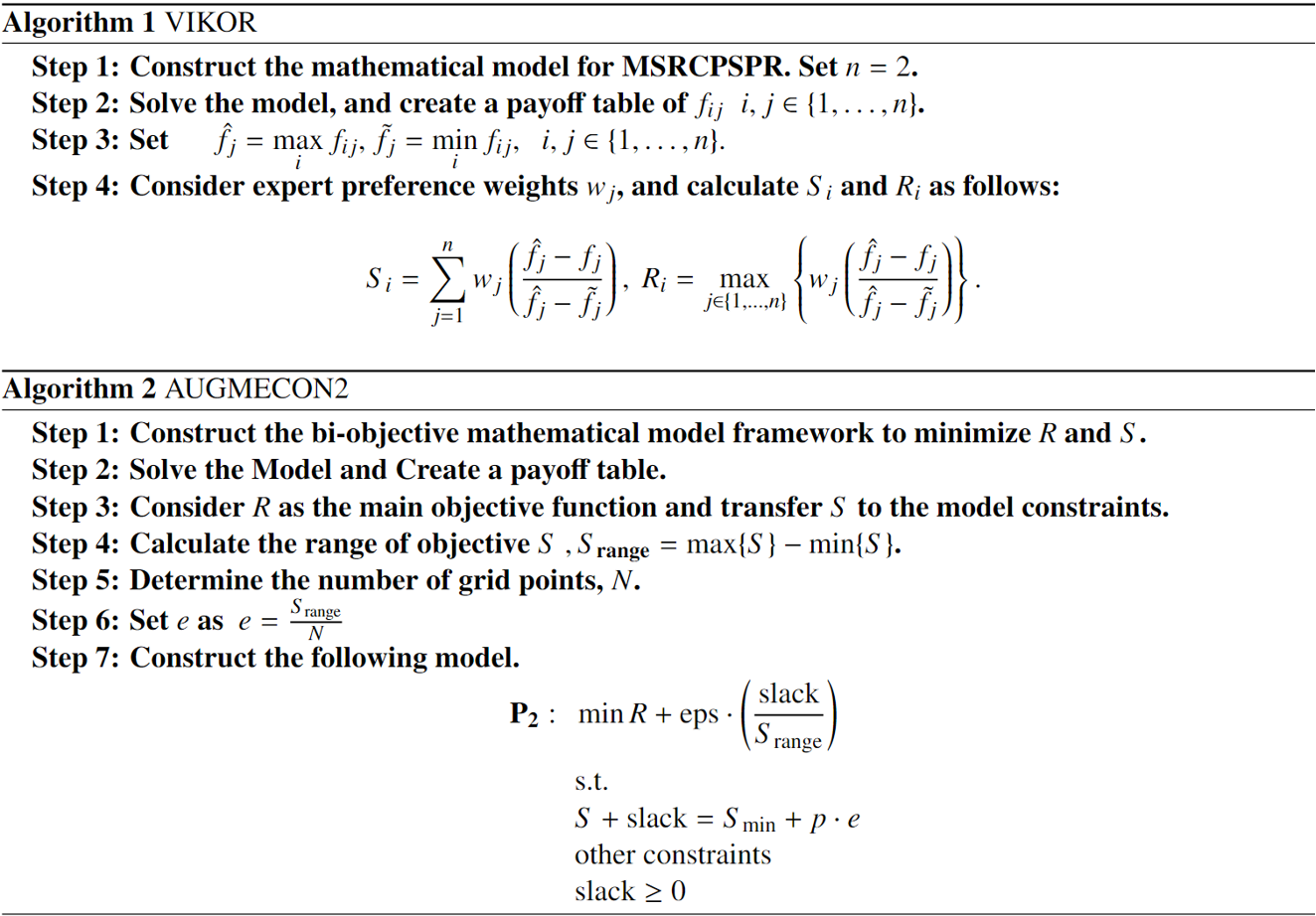}
    \caption{AUGMECON2VIKOR algorithm}
    \label{Fig1}
\end{figure}

The aim of solving the MSRCPSPR is to minimize project makespan and project costs simultaneously. Minimizing the project's makespan increases the cost, and minimizing the cost increases the makespan; therefore, the two objectives are competing. To solve the multi-objective problem, a hybrid method extended by \citet{shahabi2021routing}, called AUGMECON2VIKOR, is utilized to find Pareto solution points. This method comprises two different components, namely VIKOR method \citep{opricovic1998multicriteria} and augmented 2 epsilon-constraint  method (AUGMECON2) \citep{mavrotas2013improved}. As a selection mechanism, the VIKOR method is utilized in this paper's ranking phase of AUGMECON2. These two algorithms are demonstrated in Fig. (\ref{Fig1}).

To use the AUGMECON2VIKOR algorithm for solving the MSRCPSPR, the problem is initially solved using VIKOR algorithm to obtain $S_j$ and $R_j$ for $j \in \{1, \ldots, n\}$. Then, after setting the counter of Pareto grid points $p=1$, problem $P_2$ is constructed and solved by AUGMECON2 algorithm. The solution is added to the Pareto solution sets, and $p$ is incremented until $p \le N$. In problem $P_2$, the parameter $\text{eps} \in [10^{-6}, 10^{-3}]$  and $\text{slack}$ is the slack variable to convert the inequality constraint into an equality constraint following the notation on epsilon-constraint method.


\section{EXPERIMENTAL RESULTS AND COMPARATIVE ANALYSIS}
\label{sec:experimental results and comparative results]} 
\begin{figure}[h]
    \centering
    \includegraphics[width=0.7\linewidth]{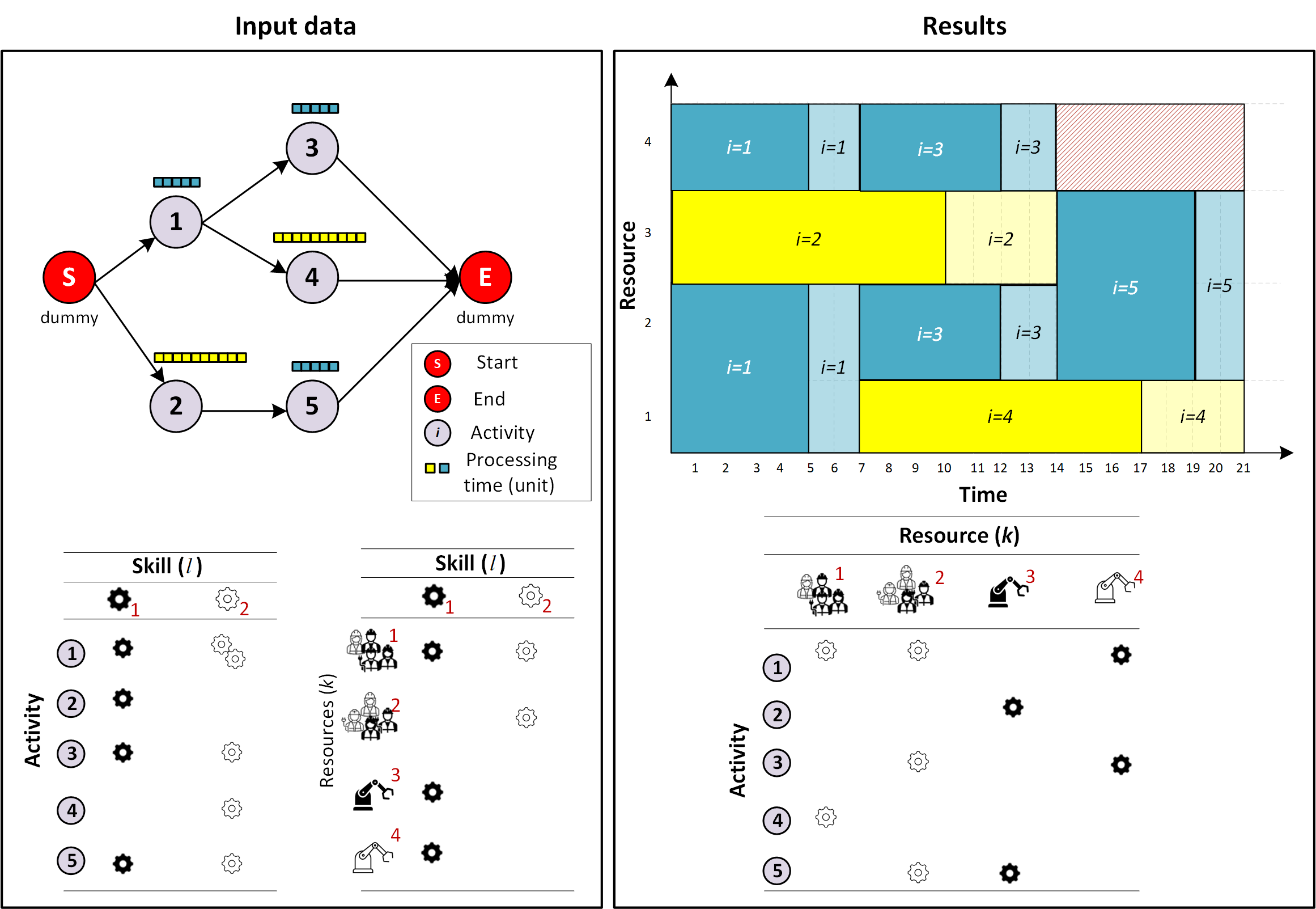}
    \caption{A numerical example}
    \label{Fig2}
\end{figure}

To validate the proposed model and solution strategy, initially, a numerical example is solved with five project activities. The representation network, input data, Gantt chart, and resource assignments are shown in Fig. (\ref{Fig2}). The first and second renewable resources are considered as human workforces and two others are machines. In addition, as shown in this figure, light color rectangles for a particular activity $i$ represent waiting times for that activity. The project is completed in 21 days. Further evaluation of this study is carried out utilizing an adapted version of problem instances j10 and j20 from the PSPLIB (\url{http://www.om-db.wi.tum.de/psplib/}), which was created by ProGen \citep{kolisch1997psplib}. Problem instances and the numerical example are solved by BARON solver in GAMS 24.1.2 interface with an 11th Gen Intel(R) Core(TM) i7-11800H @ 2.30GHz as well as 16.0 GB RAM under standard configurations.

The numerical results of solving instance j10 of the MSRCPSPR using both the AUGMECON2VIKOR algorithm and the epsilon-constraint method including the Pareto points, project makespan, and cost values are reported in Table (\ref{tab:Table1}) with fixed parameters $\upsilon_k=0.5$ and $r_k=0.5$ for both methods. The last two columns represent the percentage changes. Table (\ref{tab:Table2}) is dedicated to the sensitivity analysis of solving problem instance j10 with fixed parameter $\upsilon_k=0.5$ and two values of $r_k$. Table  (\ref{tab:Table3}) demonstrates the sensitivity analysis of solving j10 with fixed parameter $r_k=0.5$ and two values of $\upsilon_k$.


\begin{table} [H]
    \centering
    \caption{Comparative results of the objective functions for problem instance j10, with fixed parameter $\upsilon_k=0.5$ and $\upsilon=0.5$}
    \label{tab:Table1}
    \begin{tabular}{ccccccc}
        \hline
        \multirow{2}{*}{Grid points} & \multicolumn{2}{c}{AUGMECON2VIKOR} & \multicolumn{2}{c}{$\epsilon$-constraint} & \multicolumn{2}{c}{Percentage Change} \\
        \cline{2-7}
        & $Z_1$ (Makespan) & $Z_2$ (Cost) & $Z_1$ (Makespan) & $Z_2$ (Cost) & $Z_1$ (Makespan) & $Z_2$ (Cost) \\
        \hline
        1 & \textbf{48.86} & \textbf{7,080,000} & 50.06 & 7,460,000 & -2.40\% & -5.09\% \\
        2 & \textbf{52.35} & \textbf{6,740,000} & 54.87 & 6,835,000 & -4.59\% & -1.39\% \\
        3 & \textbf{56.89} & \textbf{6,240,000} & 57.42 & 6,522,500 & -0.92\% & -4.33\% \\
        4 & \textbf{63.71} & \textbf{6,190,000} & 65.69 & 6,210,000 & -3.01\% & -0.32\% \\
        5 & \textbf{72.11} & \textbf{5,760,000} & 74.25 & 5,897,500 & -2.88\% & -2.33\% \\
        6 & \textbf{73.86} & \textbf{5,540,000} & 75.11 & 5,708,000 & -1.66\% & -2.94\% \\
        7 & 74.11 & 4,960,000 & - & - & - & - \\
        \hline
    \end{tabular}
\end{table}
As shown in Table (\ref{tab:Table1}), Pareto solution points obtained from AUGMECON2VIKOR algorithm are compared to the output of the epsilon-constraint method. It is clear from the table that AUGMECON2VIKOR algorithm performs better in both finding more Pareto solution points and values of the objective functions. According to the average project makespan and cost of grid points 1 to 6, these values are reduced by approximately $2.55\%$ and $2.80\%$ in the solutions obtained from AUGMEKON2VIKOR in comparison to the solutions obtained from epsilon-constraint method. The AUGMECON2VIKOR algorithm is executed for approximately 4024 seconds to find seven Pareto solutions; however, It takes 4567 seconds for the epsilon-constraint method to find six Pareto solutions. 

\begin{figure}[H]
    \centering
    \includegraphics[width=0.584\linewidth]{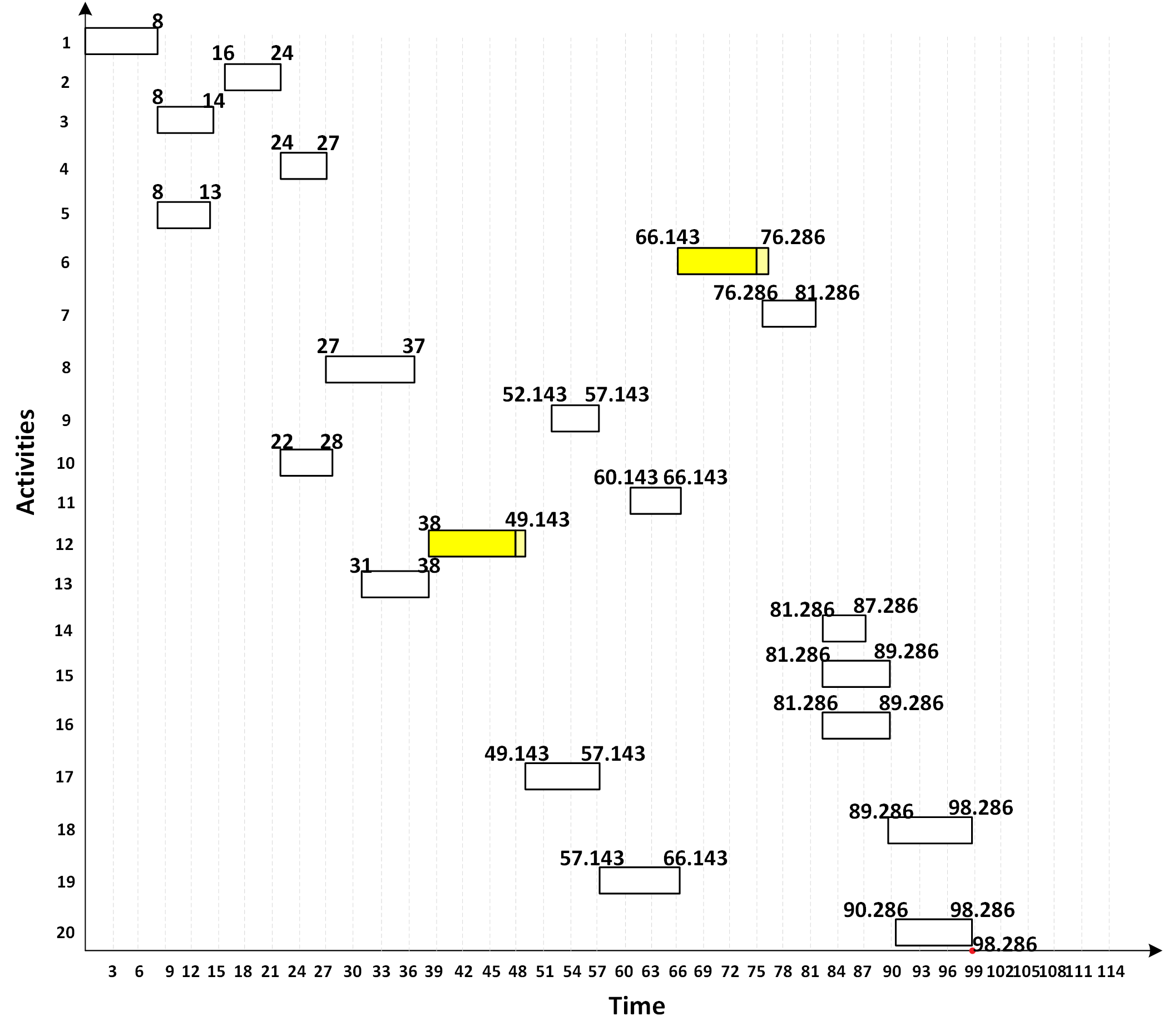} 
    \caption{Gantt chart of j20}
    \label{Fig3}
\end{figure}

To further show the capability of AUGMEKON2VIKOR algorithm, problem instance j20 is solved by this method within an execution time of nearly 6124 seconds with the outputs  $(Z_1PIS,\ Z_1NIS) =(98.29,\ 149.82)$, and $(Z_2PIS,\ Z_2NIS) =(1.78E+07,\ 1.12E+07)$, where $Z_hNIS$ and $Z_hPIS$ are the negative ideal solution (NIS) and positive ideal solution (PIS) of the $h^{th}$ objective function, respectively. It should be noted that the epsilon-constraint method is not able to find a solution under the fixed CPU time (i.e., 6500 seconds). Fig. (\ref{Fig3}) demonstrates the Gantt chart of j20 in the PIS Pareto solution point. 
In this figure, the start and finish times of each activity are represented at the beginning and end points of each box. In the Gantt chart of the provided solution, only activities 6 and 12 have waiting time which is depicted as yellow boxes in the figure.
From a sensitivity analysis point of view, it is expected that as the retrieval rate increases, the amount of retrieval time increases as well. Consequently, the length of the queue will be reduced, and waiting time will be decreased, as evidenced in Table (\ref{tab:Table2}). On the other hand, Table (\ref{tab:Table3}) certifies that an increase in the disruption rate leads to an increase in the queue length. As a result, waiting time will be increased. Furthermore, according to Table (\ref{tab:Table2}), once the retrieval rate increased by 40\%, project makespan and project cost decreased by 3.13\% and 0.25\%, respectively. In addition, the percentage change in the objective functions by increasing $\upsilon_k$ is positive in all grid points implying that the situation is worse, where project makespan increased by 2.94\% and project cost increased by 0.78\%. 

\begin{table}[H]
    \centering
    \caption{Comparative results of the objective functions for problem instance j10, with fixed parameter $\upsilon_k=0.5$ and two values of $r_k$ using the AUGMECON2VIKOR method}
    \label{tab:Table2}
    \begin{tabular}{ccccccc}
        \hline
        \multirow{2}{*}{Grid Points} & \multicolumn{2}{c}{$r_k=0.5$} & \multicolumn{2}{c}{$140\%r_k$} & \multicolumn{2}{c}{Percentage Change} \\
        \cline{2-7}
        & $Z_1$ (Makespan) & $Z_2$ (Cost) & $Z_1$ (Makespan) & $Z_2$ (Cost) & $Z_1$ (Makespan) & $Z_2$ (Cost) \\
        \hline
        1 & 48.86 & 7,080,000 & \textbf{48.6} & \textbf{7,020,000} & -0.53\% & -0.85\% \\
        2 & 52.35 & 6,740,000 & \textbf{51.81} & \textbf{6,540,000} & -1.03\% & -2.97\% \\
        3 & 56.89 & 6,240,000 & \textbf{55.26} & \textbf{6,310,000} & -2.87\% & -2.87\% \\
        4 & 63.71 & 6,190,000 & \textbf{61.6} & \textbf{6,120,000} & -3.31\% & -1.13\% \\
        5 & 72.11 & 5,760,000 & \textbf{68.59} & \textbf{5,840,000} & -4.88\% & 1.39\% \\
        6 & 73.86 & 5,540,000 & \textbf{71.21} & \textbf{5,580,000} & -3.59\% & 0.72\% \\
        7 & 74.11 & 4,960,000 & \textbf{69.89} & 4,960,000 & -5.69\% & 0.00\% \\
        \hline
    \end{tabular}
\end{table}
\begin{table}[H]
        \centering
         \caption{Comparative results of the objective functions in problem instance j10, with fixed parameter $r_k=0.5$ and two values of  $\upsilon_k$ using the AUGMECON2VIKOR method}
        \label{tab:Table3}
        \begin{tabular}{ccccccc}
        \hline
        \multirow{2}{*}{\textbf{Grid Points}} & \multicolumn{2}{c}{\textbf{$\upsilon_k=0.5$}} & \multicolumn{2}{c}{\textbf{$140\%\upsilon_k$}} & \multicolumn{2}{c}{Percentage Change} \\
        \cline{2-7}
        & $Z_1$ (Makespan) &$Z_2$ (Cost) & $Z_1$ (Makespan) & $Z_2$ (Cost) & $Z_1$ (Makespan) & $Z_2$ (Cost) \\
        \hline
        1 & \textbf{48.86} & \textbf{6,740,000} & 54.47 & 6,770,000 & 4.05\% & 0.45\% \\
        3 & \textbf{56.89} & \textbf{6,240,000} & 58.81 & 6,440,000 & 3.37\% & 3.21\% \\
        4 & \textbf{63.71} & \textbf{6,190,000} & 65.85 & 6,230,000 & 3.36\% & 0.65\% \\
        5 & \textbf{72.11} & \textbf{5,760,000} & 73.52 & 5,790,000 & 1.96\% & 0.52\% \\
        6 & \textbf{73.86} & \textbf{5,540,000} & 75.71 & 5,560,000 & 2.50\% & 0.36\% \\
        7 & \textbf{74.11} & 4,960,000 & 77.77 & 4,960,000 & 4.94\% & 0.00\% \\
        \hline
    \end{tabular}
\end{table}
\sethlcolor{yellow}
\section{CONCLUSION AND FURTHER RESEARCH}
\label{sec:experimental results} 
This paper introduces a novel bi-objective MINLP model for the MSRCPSP that aims to minimize project makespan and cost while considering reliability considerations. The model addresses a gap in the literature by treating renewable resources, including human workforce and machines, as reliability factors within the MSRCPSP. To take into consideration disruptions in complex manufacturing projects, the $M/M/1$ queueing system is integrated to represent the reliability aspect of the model. The model's effectiveness is validated through a numerical example with five tasks. The results compare the performance of AUGMECON2VIKOR algorithm to the widely used epsilon-constraint method solving adapted problem instances j10 and j20 of the MSRCPSPR. A comparative analysis shows a positive correlation between increasing the repair rate of machines and objective functions (i.e. makespan and cost). On the other hand, increasing the disruption rate of machines leads to a rise in both project makespan (2.94\%) and project cost (0.78\%). These findings highlight the significant impact of incorporating reliability constraints on project completion time and cost optimization. This study focuses on scheduling a single project. In today's competitive markets, companies require multi-project scheduling. Such scenarios involve higher complexity and potentially greater disruption rates. In addition, the study assumes each activity has one execution mode. However, allowing for multiple execution modes is an important aspect of the MSRCPSPR for representing more realistic scenarios. Investigating large-scale problems using meta-heuristic algorithms and multi-project multi-execution mode cases is an important extension of the current work.
\bibliographystyle{elsarticle-num-names} 
\bibliography{cie-refs}





\end{document}